\documentclass[12pt]{article}
\usepackage{amssymb}
\usepackage{amsmath}

\input{tcilatex}

\begin{document}

\begin{center}
\textbf{A characterization of compact operators via the non-connectedness of
the attractors of a family of IFSs}

\bigskip

by

\bigskip

ALEXANDRU\ MIHAIL (Bucharest) and

RADU\ MICULESCU (Bucharest)

\bigskip
\end{center}

{\small Abstract}. {\small In this paper we present a result which
establishes a connection between the theory of compact operators and the
theory of iterated function systems. For a Banach space }$X${\small , }$S$%
{\small \ and }$T${\small \ bounded linear operators from }$X${\small \ to }$%
X${\small \ such that }$\left\Vert S\right\Vert ,\left\Vert T\right\Vert <1$%
{\small \ and }$w\in X${\small , let us consider the IFS }$%
S_{w}=(X,f_{1},f_{2})${\small , where }$f_{1},f_{2}:X\rightarrow X${\small \
are given by }$f_{1}(x)=S(x)${\small \ and }$f_{2}(x)=T(x)+w${\small , for
all }$x\in X${\small . On one hand we prove that if the operator }$S${\small %
\ is compact, then there exists a family }$(K_{n})_{n\in \mathbb{N}}${\small %
\ of compact subsets of }$X${\small \ such that }$A_{\mathcal{S}_{w}}$%
{\small \ is\ not connected, for all }$w\in H-\underset{n\in \mathbb{N}}{%
\cup }K_{n}${\small . One the other hand we prove that if }$H$ {\small is an
infinite dimensional Hilbert space, then a bounded linear operator }$%
S:H\rightarrow H${\small \ having the property that }$\left\Vert
S\right\Vert <1${\small \ is compact provided that for every bounded linear
operator }$T:H\rightarrow H${\small \ such that }$\left\Vert T\right\Vert <1$%
{\small \ there exists a sequence }$(K_{T,n})_{n}${\small \ of compact
subsets of }$H${\small \ such that }$A_{\mathcal{S}_{w}}${\small \ is\ not
connected for all }$w\in H-\underset{n}{\cup }K_{T,n}${\small .
Consequently, given an infinite dimensional Hilbert space }$H$,{\small \
there exists a complete characterization of the compactness of an operator }$%
S:H\rightarrow H${\small \ by means of the non-connectedness of the
attractors of a family of IFSs related to the given operator.}

\bigskip

\textbf{1. Introduction. }IFSs were introduced in their present form by John
Hutchinson (see [9]) and popularized by Michael Barnsley (see [2]). They are
one of the most common and most general ways to generate fractals. Although
the fractals sets are defined by means of measure theory concepts (see [7]),
they have very interesting topological properties. The connectivity of the
attractor of an iterated function system has been studied, for example, in
[14] (for the case of an iterated multifunction system) and in [6] (for the
case of an infinite iterated function system).

It is well known the role of the compact operators theory in functional

---------------------------

\textit{2000 Mathematical\ Subject Classification}{\small : }Primary 28A80,
47B07; Secondary 54D05


\bigskip \newpage analysis and, in particular, in the theory of the integral
equations. In this frame, a natural question is to provide equivalent
characterizations for compact operators. Let us mention some results on this
direction. A bounded operator $T$ on a separable Hilbert space $H$ is
compact if and only if $\underset{n\rightarrow \infty }{\lim }%
<Te_{n},e_{n}>=0$ (or equivalently $\underset{n\rightarrow \infty }{\lim }%
Te_{n}=0$), for each orthonormal basis $\{e_{n}\}$ for $H$ (see [1], [8],
[16] and [17]) if and only if every orthonormal basis $\{e_{n}\}$ for $H$
has a rearrangement $\{e_{\sigma (n)}\}$ such that $\sum \frac{1}{n}%
\left\Vert Te_{\sigma (n)}\right\Vert <\infty $ (see [18]). In a more
general framework, in [10] a characterization of the compact operators on a
fixed Banach space in terms of a construction due to J.J.M. Chadwick and
A.W. Wickstead (see [3]) is presented and in [11] a purely structural
characterization of compact elements in a $C^{\ast }$ algebra is given.

In contrast to the above mentioned characterizations of the compact
operators which are confined to the framework of the functional analysis, in
this paper we present such a characterization by means of the
non-connectedness of the attractors of a family of IFSs related to the
considered operator.

In this way we establish an unexpected connection between the theory of
compact operators and the theory of iterated function systems.

\bigskip

\textbf{2.} \textbf{Preliminary results. }In this paper, for a function $f$\
and $n\in \mathbb{N}$, by $f^{[n]}$ we mean the composition of $f$ by itself 
$n$ times.

\bigskip

DEFINITION\ 2.1. Let $(X,d)$\ be a metric space. A function $f:X\rightarrow
X $\ is called a \textit{contraction} in case there exists $\lambda \in
(0,1) $ such that%
\begin{equation*}
d(f(x),f(y))\leq \lambda d(x,y)\text{,}
\end{equation*}%
for all $x,y\in X$.

\bigskip

THEOREM\ 2.2 (The Banach-Cacciopoli-Picard contraction principle). \textit{%
If }$X$\textit{\ is a complete metric space, then for each contraction }$%
f:X\rightarrow X$\textit{\ there exists a unique fixed point }$x^{\ast }$%
\textit{\ of }$f$\textit{.}

\textit{Moreover }%
\begin{equation*}
x^{\ast }=\underset{n\rightarrow \infty }{\lim }f^{[n]}(x_{0})\text{\textit{,%
}}
\end{equation*}%
\textit{\ for each }$x_{0}\in X$\textit{.}

\bigskip

NOTATION. Given a metric space $(X,d)$, by $K(X)$\ we denote the set of
non-empty compact subsets of $X$.

\bigskip

DEFINITION\ 2.3. For a metric space $(X,d)$, the function $h:K(X)\times
K(X)\rightarrow \lbrack 0,+\infty )$ defined by

\begin{equation*}
h(A,B)=\max (d(A,B),d(B,A))=
\end{equation*}

\begin{equation*}
=\inf \{r\in \lbrack 0,\infty ):A\subseteq B(B,r)\text{ and }B\subseteq
B(A,r)\}\text{,}
\end{equation*}%
where%
\begin{equation*}
B(A,r)=\{x\in X:d(x,A)<r\}
\end{equation*}%
and

\begin{equation*}
d(A,B)=\underset{x\in A}{\sup }d(x,B)=\underset{x\in A}{\sup }(\underset{%
y\in B}{\inf }d(x,y))\text{,}
\end{equation*}%
turns out to be a metric which is called the \textit{Hausdorff-Pompeiu metric%
}.

\bigskip

REMARK 2.4. The metric space $(K(X),h)$ is complete, provided that $(X,d)$
is a complete metric space.

\bigskip

DEFINITION\ 2.5. Let $(X,d)$\ be a complete metric space. An \textit{%
iterated function system} (for short an IFS) on $X$, denoted by $%
S=(X,(f_{k})_{k\in \{1,2,...,n\}})$, consists of a finite family of
contractions $(f_{k})_{k\in \{1,2,...,n\}}$, $f_{k}:X\rightarrow X$.

\bigskip

THEOREM\ 2.6. \textit{Given} $\mathcal{S}=(X,(f_{k})_{k\in \{1,2,...,n\}})$ 
\textit{an iterated function system on }$X$, \textit{the function} $F_{%
\mathcal{S}}:K(X)\rightarrow K(X)$ \textit{defined by}

\begin{equation*}
F_{\mathcal{S}}(C)=\underset{k=1}{\overset{n}{\cup }}f_{k}(C)\text{,}
\end{equation*}%
\textit{for all} $C\in K(X)$, \textit{which is called the set function
associated to }$\mathcal{S}$\textit{, turns out to be a contraction and its
unique fixed point, denoted by }$A_{\mathcal{S}}$\textit{, is called the
attractor of the IFS }$\mathcal{S}$.

\bigskip

REMARK 2.7. For each $i\in \{1,2,...,n\}$, the fixed point of $f_{i}$\ is an
element of $A_{\mathcal{S}}$.

\bigskip

REMARK 2.8. If $A\in K(X)$ has the property that $F_{\mathcal{S}%
}(A)\subseteq A$, then $A_{\mathcal{S}}\subseteq A$.

\bigskip

$\mathtt{Proof}$. The proof is similar to the one of Lemma 3.6 from [13]. $%
\square $

\bigskip

DEFINITION\ 2.9. Let $(X,d)$\ be a metric space and $(A_{i})_{i\in I}$\ a
family of nonempty subsets of $X$. The \textit{family} $(A_{i})_{i\in I}$\
is said to be \textit{connected} if for every $i,j\in I$, there exist $n\in 
\mathbb{N}$\ and $\{i_{1},i_{2},...,i_{n}\}\subseteq I$\ such that $i_{1}=i$%
, $i_{n}=j$\ and $A_{i_{k}}\cap A_{i_{k+1}}\neq \emptyset $\ for every $k\in
\{1,2,..,n-1\}$.

\bigskip

THEOREM\ 2.10\textit{\ }(see [12], Theorem 1.6.2, page 33)\textit{. Given an
IFS }$\mathcal{S}=(X,(f_{k})_{k\in \{1,2,...,n\}})$\textit{, where }$(X,d)$%
\textit{\ is a complete metric space, the following statements are
equivalent:}

\textit{1)\ the family }$(f_{i}(A_{\mathcal{S}}))_{i\in \{1,2,...,n\}}$%
\textit{\ is connected;}

\textit{2) }$A_{\mathcal{S}}$\textit{\ is arcwise connected.}

\textit{3) }$A_{\mathcal{S}}$\textit{\ is connected.}

\bigskip

PROPOSITION 2.11. \textit{For a given complete metric space }$(X,d)$\textit{%
, let us consider the IFSs }$\mathcal{S}=(X,f_{1},f_{2})$\textit{\ and }$%
\mathcal{S}^{^{\prime }}=(X,f_{1}^{[m]},f_{2})$\textit{, where }$m\in 
\mathbb{N}$\textit{.}

\textit{If }$A_{\mathcal{S}^{^{\prime }}}$\textit{\ is connected, then }$A_{%
\mathcal{S}}$\textit{\ is connected.}

\bigskip

$\mathtt{Proof}$. Since $F_{\mathcal{S}^{^{\prime }}}(A_{\mathcal{S}%
})=f_{1}^{[m]}(A_{\mathcal{S}})\cup f_{2}(A_{\mathcal{S}})\subseteq A_{%
\mathcal{S}}$, we get (using Remark 2.8) $A_{\mathcal{S}^{^{\prime
}}}\subseteq A_{\mathcal{S}}$ and hence $f_{2}(A_{\mathcal{S}^{^{\prime
}}})\subseteq f_{2}(A_{\mathcal{S}})$. Because $f_{1}^{[m]}(A_{\mathcal{S}%
^{^{\prime }}})\subseteq f_{1}(A_{\mathcal{S}})$, it follows that $%
f_{1}^{[m]}(A_{\mathcal{S}^{^{\prime }}})\cap f_{2}(A_{\mathcal{S}^{^{\prime
}}})\subseteq f_{1}(A_{\mathcal{S}})\cap f_{2}(A_{\mathcal{S}})$ (*). Since $%
A_{\mathcal{S}^{^{\prime }}}$\textit{\ }is connected, taking into account
Theorem 2.10, we deduce that $f_{1}^{[m]}(A_{\mathcal{S}^{^{\prime }}})\cap
f_{2}(A_{\mathcal{S}^{^{\prime }}})\neq \varnothing $, which, using $(\ast )$%
, implies that $f_{1}(A_{\mathcal{S}})\cap f_{2}(A_{\mathcal{S}})\neq
\varnothing $. Then, using again Theorem 2.10, we infer that $A_{\mathcal{S}%
} $\textit{\ }is connected. $\square $

\bigskip

PROPOSITION 2.12\textbf{\ }(see [5], page 238, lines 11-12). \textit{Assume
that }$H$\textit{\ is a Hilbert space. Let us consider a self-adjoint
operator }$N:H\rightarrow H$\textit{\ and }$E$ \textit{its spectral
decomposition. Then for each }$\lambda \in \mathbb{R}$\textit{\ we have}%
\begin{equation*}
NE((-\infty ,\lambda ))\leq \lambda E((-\infty ,\lambda ))
\end{equation*}%
\textit{and}%
\begin{equation*}
\lambda E((\lambda ,\infty ))\leq NE((\lambda ,\infty ))\text{,}
\end{equation*}%
\textit{for all }$\lambda \in \mathbb{R}$\textit{.}

\bigskip

PROPOSITION 2.13\textbf{\ }(see [5], page 226, Observation 7). \textit{%
Assume that }$H$\textit{\ is a Hilbert space. Let us consider two
self-adjoint operators }$N_{1},N_{2}:H\rightarrow H$\textit{.}

\textit{If}%
\begin{equation*}
0\leq N_{1}\leq N_{2}\text{,}
\end{equation*}%
\textit{then}%
\begin{equation*}
\left\Vert N_{1}\right\Vert \leq \left\Vert N_{2}\right\Vert \text{.}
\end{equation*}

\bigskip

PROPOSITION 2.14\textbf{\ }(see [19], ex. 25, page 344). \textit{Assume that 
}$H$\textit{\ is a Hilbert space. Let us consider a normal operator }$%
N:H\rightarrow H$\textit{, }$g$ \textit{a bounded Borel function on }$\sigma
(N)$\textit{\ and }$S=g(T)$\textit{. If }$E_{N}$\textit{\ and }$E_{S}$%
\textit{\ are the spectral decomposition of }$N$\textit{\ and }$S$\textit{,
then }%
\begin{equation*}
E_{S}(\omega )=E_{N}(g^{-1}(\omega ))\text{,}
\end{equation*}%
\textit{for every Borel set }$\omega \subseteq \sigma (S)$\textit{.}

\bigskip

PROPOSITION 2.15 (see [4], Proposition 4.1, page 278). \textit{Assume that }$%
H$\textit{\ is a Hilbert space. Let us consider a normal operator} $%
N:H\rightarrow H$ \textit{and} $E$ \textit{its spectral decomposition.\ Then 
}$N$\textit{\ is compact if and only if }$E(\{z\mid \left\vert z\right\vert
>\varepsilon \})$\textit{\ has finite rank, for every }$\varepsilon >0$%
\textit{.}

\bigskip

PROPOSITION 2.16. \textit{Assume that }$H$\textit{\ is a Hilbert space. Let
us consider a bounded linear operator }$A:H\rightarrow H$ \textit{which is
invertible. Then }$Id_{H}-A^{\ast }A$\textit{\ is compact if and only if }$%
Id_{H}-AA^{\ast }$\textit{\ is compact.}

\bigskip

$\mathtt{Proof}$. According to the well known polar decomposition theorem
there exists an unitary operator $U:H\rightarrow H$ and a positive operator $%
P:H\rightarrow H$ such that $P^{2}=A^{\ast }A$ and $A=UP$. Then%
\begin{equation*}
Id_{H}-AA^{\ast }=Id_{H}-UP(UP)^{\ast }=Id_{H}-UPP^{\ast }U^{\ast
}=Id_{H}-UP^{2}U^{\ast }=
\end{equation*}%
\begin{equation*}
=UU^{\ast }-UP^{2}U^{\ast }=U(Id_{H}-P^{2})U^{\ast }=U(Id_{H}-A^{\ast
}A)U^{\ast }\text{.}
\end{equation*}

Hence $Id_{H}-AA^{\ast }=U(Id_{H}-A^{\ast }A)U^{\ast }$ and $Id_{H}-A^{\ast
}A=U^{\ast }(Id_{H}-AA^{\ast })U$. From the last two relations we obtain the
conclusion. $\square $

\bigskip

COROLLARY 2.17. \textit{Assume that }$H$\textit{\ is a Hilbert space. Let us
consider a bounded linear operator }$S:H\rightarrow H$ \textit{such that }$%
\left\Vert S\right\Vert <1$\textit{. Then }$S+S^{\ast }-SS^{\ast }$\textit{\
is compact if and only if }$S+S^{\ast }-S^{\ast }S$\textit{\ is compact.}

\bigskip

$\mathtt{Proof}$. The operator $A=Id_{H}-S$ is invertible since $\left\Vert
S\right\Vert <1$. According to Proposition 2.16 $Id_{H}-A^{\ast }A$\textit{\ 
}is compact if and only if $Id_{H}-AA^{\ast }$ is compact i.e. $S+S^{\ast
}-SS^{\ast }$ is compact if and only if $S+S^{\ast }-S^{\ast }S$ is compact. 
$\square $

\bigskip

PROPOSITION 2.18\textbf{\ }(see [19], ex. 14, page 324). \textit{Assume that 
}$H$\textit{\ is a Hilbert space and let us consider a bounded linear
operator }$S:H\rightarrow H$\textit{. If }$S^{\ast }S$\textit{\ is a compact
operator, then }$S$\textit{\ is compact.}

\bigskip

\textbf{3.}\ \textbf{A sufficient condition for the compactness of an
operator. }In this section, $H$ is an infinite-dimensional Hilbert space. We
shall use the notation $Id_{H}$ for the function $Id_{H}:H\rightarrow H$,
given by $Id_{H}(x)=x$, for all $x\in H$ . If $S$ and $T$ are bounded linear
operators from $H$ to $H$ such that $\left\Vert S\right\Vert ,\left\Vert
T\right\Vert <1$, then $S$ and $T$ are contractions. For $w\in X$, we
consider the IFS $S_{w}=(X,f_{1},f_{2})$, where $f_{1},f_{2}:X\rightarrow X$
are given by $f_{1}(x)=S(x)$ and $f_{2}(x)=T(x)+w$, for all $x\in X$.

\bigskip

THEOREM\ 3.1. \textit{In the preceding framework, let us consider a bounded
linear operator }$S:H\rightarrow H$ \textit{satisfying the condition} $%
\left\Vert S\right\Vert <1$\textit{. If for every bounded linear operator }$%
T:H\rightarrow H$\textit{\ such that }$\left\Vert T\right\Vert <1$\textit{\
there exists a sequence }$(K_{T,n})_{n}$\textit{\ of compact subsets of }$H$%
\textit{\ having the property that }$A_{\mathcal{S}_{w}}$\textit{\ is\ not
connected for all }$w\in H-\underset{n}{\cup }K_{T,n}$\textit{, then the
operator }$S$\textit{\ is compact.}

\bigskip

$\mathtt{Proof}$. For each $m\in \mathbb{N}$ let us consider the bounded
linear operator $U=S^{[m]}$. Obviously $\left\Vert U\right\Vert <1$. Let us
consider $P_{\varepsilon }=E((-\infty ,1-\varepsilon ))$ and $\overset{\sim }%
{P_{\varepsilon }}=E((1+\varepsilon ,\infty ))$, where $E$ is the spectral
decomposition of the positive (so self-adjoint, so normal) bounded linear
operator 
\begin{equation*}
N=(Id_{H}-U)^{\ast }(Id_{H}-U)=Id_{H}-U-U^{\ast }+U^{\ast }U\text{.}
\end{equation*}

We claim\textit{\ that }$P_{\varepsilon }$\textit{\ has finite rank for
every }$\varepsilon >0$.

Indeed, if there is to be an $\varepsilon _{0}>0$\textit{\ }such that $%
P_{\varepsilon _{0}}$ has infinite rank, then\textit{\ }let us consider the
operator $T=(Id_{H}-U)P_{\varepsilon _{0}}$ and remark that%
\begin{equation*}
NP_{\varepsilon _{0}}=NP_{\varepsilon _{0}}^{2}=NP_{\varepsilon _{0}}^{\ast
}P_{\varepsilon _{0}}=P_{\varepsilon _{0}}^{\ast }NP_{\varepsilon
_{0}}=P_{\varepsilon _{0}}^{\ast }((Id_{H}-U)^{\ast
}(Id_{H}-U))P_{\varepsilon _{0}}=
\end{equation*}%
\begin{equation*}
=((Id_{H}-U)P_{\varepsilon _{0}})^{\ast }((Id_{H}-U)P_{\varepsilon
_{0}})\geq 0\text{.}
\end{equation*}%
Hence, according to Proposition 2.12, we have $0\leq NP_{\varepsilon
_{0}}\leq (1-\varepsilon _{0})P_{\varepsilon _{0}}$ and therefore, using
Proposition 2.13, it follows that $\left\Vert NP_{\varepsilon
_{0}}\right\Vert \leq 1-\varepsilon _{0}$. Consequently we obtain 
\begin{equation*}
\left\Vert T\right\Vert ^{2}=\left\Vert T^{\ast }T\right\Vert =\left\Vert
(Id_{H}-U)P_{\varepsilon _{0}})^{\ast }(Id_{H}-U)P_{\varepsilon
_{0}}\right\Vert =
\end{equation*}%
\begin{equation*}
=\left\Vert P_{\varepsilon _{0}}^{\ast }(Id_{H}-U)^{\ast
}(Id_{H}-U)P_{\varepsilon _{0}}\right\Vert =\left\Vert P_{\varepsilon
_{0}}NP_{\varepsilon _{0}}\right\Vert \leq
\end{equation*}%
\begin{equation*}
\leq \left\Vert P_{\varepsilon _{0}}\right\Vert \left\Vert NP_{\varepsilon
_{0}}\right\Vert =\left\Vert NP_{\varepsilon _{0}}\right\Vert \leq
1-\varepsilon _{0}
\end{equation*}%
and thus%
\begin{equation*}
\left\Vert T\right\Vert \leq \sqrt{1-\varepsilon _{0}}<1\text{.}
\end{equation*}

For $w\in H$, let us consider, besides $\mathcal{S}_{w}$, the IFS $\mathcal{S%
}_{w}^{^{\prime }}=(H,f,f_{2})$, where $f:H\rightarrow H$ is given by $%
f(x)=U(x)$, for all $x\in H$.

Now let us choose an arbitrary $w\in (Id_{H}-T)P_{\varepsilon _{0}}(H)$. On
one hand, since $0$ is the fixed point of $f$, using Remark 2.7, we infer
that $0\in A_{\mathcal{S}_{w}}$. On the other hand, using the same argument,
we get that $e$, the fixed point of $f_{2}$, belongs to $A_{\mathcal{S}_{w}}$%
, that is $e=U^{-1}(w)=(Id_{H}-T)^{-1}(w)\in A_{\mathcal{S}_{w}^{^{\prime
}}} $. Since $f(e)=f_{2}(0)=w$, we obtain $w\in f(A_{\mathcal{S}%
_{w}^{^{\prime }}})\cap f_{2}(A_{\mathcal{S}_{w}^{^{\prime }}})$, which
implies $f(A_{\mathcal{S}_{w}^{^{\prime }}})\cap f_{2}(A_{\mathcal{S}%
_{w}^{^{\prime }}})\neq \emptyset $, and therefore, according to Theorem
2.10, $A_{\mathcal{S}_{w}^{^{\prime }}}$is connected. We conclude (using
Proposition 2.11) that $A_{\mathcal{S}_{w}}$ is connected.

Consequently there exists a bounded linear operator $T:H\rightarrow H$
having $\left\Vert T\right\Vert <1$ such that $A_{\mathcal{S}_{w}}$ is
connected for every $w\in (Id_{H}-T)P_{\varepsilon _{0}}(H)$.

According to the hypothesis there exists a sequence $(K_{T,n})_{n}$\textit{\ 
}of compact subsets of $H$\ having the property that $A_{\mathcal{S}_{w}}$\
is\ not connected, for all $w\in H-\underset{n}{\cup }K_{T,n}m.$

Therefore we obtain $(Id_{H}-T)P_{\varepsilon _{0}}(H)\subseteq \underset{n}{%
\cup }K_{T,n}$ which (taking into account the fact that $(Id_{H}-T)P_{%
\varepsilon _{0}}(H)$ is infinite dimensional, that the closed unit ball in
a normed linear space $X$ is compact if and only if $X$ is infinite
dimensional and Baire's theorem) generates a contradiction.

We assert \textit{that }$\overset{\sim }{P_{\varepsilon }}$\textit{\ has
finite rank for every }$\varepsilon >0$.

Indeed, if by contrary we suppose that there exists $\varepsilon _{0}>0$%
\textit{\ }such that $\overset{\sim }{P_{\varepsilon _{0}}}$ has infinite
rank, let $R_{\varepsilon _{0}}$ designates the orthogonal projection of $H$
onto $(Id_{H}-U)\overset{\sim }{P_{\varepsilon _{0}}}(H)$ and let us
consider the bounded linear operator $T=(Id_{H}-U)^{-1}R_{\varepsilon _{0}}$%
. Based upon Proposition 2.12, we have 
\begin{equation*}
N\overset{\sim }{P_{\varepsilon _{0}}}=(Id_{H}-U)^{\ast }(Id_{H}-U)\overset{%
\sim }{P_{\varepsilon _{0}}}\geq (1+\varepsilon _{0})\overset{\sim }{%
P_{\varepsilon _{0}}}\text{,}
\end{equation*}%
which implies that 
\begin{equation*}
\left\Vert (Id_{H}-U)\overset{\sim }{P_{\varepsilon _{0}}}(x)\right\Vert
^{2}=<N\overset{\sim }{P_{\varepsilon _{0}}}(x),\overset{\sim }{%
P_{\varepsilon _{0}}}(x)>\geq (1+\varepsilon _{0})\left\Vert \overset{\sim }{%
P_{\varepsilon _{0}}}(x)\right\Vert ^{2}\text{,}
\end{equation*}%
i.e.%
\begin{equation}
\sqrt{1+\varepsilon _{0}}\left\Vert \overset{\sim }{P_{\varepsilon _{0}}}%
(x)\right\Vert \leq \left\Vert (Id_{H}-U)\overset{\sim }{P_{\varepsilon _{0}}%
}(x)\right\Vert \text{,}  \tag{0}
\end{equation}%
for each $x\in H$. So, as for each $u\in H$ there exists $x_{u}\in H$ such
that $R_{\varepsilon _{0}}(u)=(Id_{H}-U)\overset{\sim }{P_{\varepsilon _{0}}}%
(x_{u})$, we infer that%
\begin{equation*}
\left\Vert T(u)\right\Vert =\left\Vert (Id_{H}-U)^{-1}R_{\varepsilon
_{0}}(u)\right\Vert =\left\Vert (Id_{H}-U)^{-1}(Id_{H}-U)\overset{\sim }{%
P_{\varepsilon _{0}}}(x_{u})\right\Vert =
\end{equation*}%
\begin{equation*}
=\left\Vert \overset{\sim }{P_{\varepsilon _{0}}}(x_{u})\right\Vert \overset{%
(0)}{\leq }\frac{1}{\sqrt{1+\varepsilon _{0}}}\left\Vert (Id_{H}-U)\overset{%
\sim }{P_{\varepsilon _{0}}}(x)\right\Vert =
\end{equation*}%
\begin{equation*}
=\frac{1}{\sqrt{1+\varepsilon _{0}}}\left\Vert R_{\varepsilon
_{0}}(u)\right\Vert \leq \frac{1}{\sqrt{1+\varepsilon _{0}}}\left\Vert
R_{\varepsilon _{0}}\right\Vert \left\Vert u\right\Vert =\frac{1}{\sqrt{%
1+\varepsilon _{0}}}\left\Vert u\right\Vert
\end{equation*}%
i.e. $\left\Vert T(u)\right\Vert \leq \frac{1}{\sqrt{1+\varepsilon _{0}}}%
\left\Vert u\right\Vert $, for each $u\in H$, which takes on the form%
\begin{equation*}
\left\Vert T\right\Vert \leq \frac{1}{\sqrt{1+\varepsilon _{0}}}<1\text{.}
\end{equation*}

For $w\in H$, let us consider, besides $\mathcal{S}_{w}$, the IFS $\mathcal{S%
}_{w}^{^{\prime }}=(H,f,f_{2})$, where $f:H\rightarrow H$ is given by $%
f(x)=U(x)$, for all $x\in H$.

Now let us choose an arbitrary $w\in (Id_{H}-T)\overset{\sim }{%
P_{\varepsilon _{0}}}(H)$. Then there exists $u\in H$ such that $w=(Id_{H}-T)%
\overset{\sim }{P_{\varepsilon _{0}}}(u)$. On one hand, since $0$ is the
fixed point of $f$, using Remark 2.7, we infer that $0\in A_{\mathcal{S}%
_{w}^{^{\prime }}}$. On the other hand, using the same argument, we get that 
$e$ (the fixed point of $f_{2}$) belongs to $A_{\mathcal{S}_{w}^{^{\prime
}}} $, that is $e=U^{-1}(w)=(Id_{H}-T)^{-1}(w)\in A_{\mathcal{S}%
_{w}^{^{\prime }}}$, and therefore $f(e)\in A_{\mathcal{S}_{w}^{^{\prime }}}$%
. Since $f(0)=0 $, on one hand we infer that%
\begin{equation}
0\in f(A_{\mathcal{S}_{w}^{^{\prime }}})\text{.}  \tag{1}
\end{equation}

On the other hand we have%
\begin{equation*}
f_{2}(f(e))=TU(e)+w=TU(Id_{H}-T)^{-1}(w)+(Id_{H}-T)(Id_{H}-T)^{-1}(w)=
\end{equation*}%
\begin{equation*}
=(Id_{H}-T(Id_{H}-U))(Id_{H}-T)^{-1}(w)=
\end{equation*}%
\begin{equation*}
=(Id_{H}-T(Id_{H}-U))(Id_{H}-T)^{-1}(Id_{H}-T)\overset{\sim }{P_{\varepsilon
_{0}}}(u)=
\end{equation*}%
\begin{equation*}
=(Id_{H}-T(Id_{H}-U))\overset{\sim }{P_{\varepsilon _{0}}}(u)=\overset{\sim }%
{P_{\varepsilon _{0}}}(u)-(Id_{H}-U)^{-1}R_{\varepsilon _{0}}(Id_{H}-U))%
\overset{\sim }{P_{\varepsilon _{0}}}(u)=
\end{equation*}%
\begin{equation*}
=\overset{\sim }{P_{\varepsilon _{0}}}(u)-(Id_{H}-U)^{-1}(Id_{H}-U))\overset{%
\sim }{P_{\varepsilon _{0}}}(u)=0\text{,}
\end{equation*}%
so%
\begin{equation}
0\in f_{2}(A_{\mathcal{S}_{w}^{^{\prime }}})\text{.}  \tag{2}
\end{equation}

From $(1)$ and $(2)$ we obtain $0\in f(A_{\mathcal{S}_{w}^{^{\prime }}})\cap
f_{2}(A_{\mathcal{S}_{w}^{^{\prime }}})$, i.e. $f(A_{\mathcal{S}%
_{w}^{^{\prime }}})\cap f_{2}(A_{\mathcal{S}_{w}^{^{\prime }}})\neq
\varnothing $, so, relying on Theorem 2.10, $A_{S_{w}^{^{\prime }}}$ is
connected. We appeal to Proposition 2.11 to deduce that $A_{\mathcal{S}_{w}}$%
is connected.

Consequently there exists a bounded linear operator $T:H\rightarrow H$
having $\left\Vert T\right\Vert <1$ such that $A_{\mathcal{S}_{w}}$is
connected for every $w\in (Id_{H}-T)\overset{\sim }{P_{\varepsilon _{0}}}(H)$%
.

Taking into account the hypothesis there exists a sequence $(K_{T,n})_{n}$%
\textit{\ }of compact subsets of $H$\ having the property that $A_{S_{w}}$\
is\ not connected for all $w\in H-\underset{n}{\cup }K_{T,n}m.$

Thus we obtain the inclusion $(Id_{H}-T)\overset{\sim }{P_{\varepsilon _{0}}}%
(H)\subseteq \underset{n}{\cup }K_{T,n}$ which generates a contradiction by
invoking the same arguments that we used in the final part of the previous
claim's proof.

Now\textit{\ we state that }$Id_{H}-(Id_{H}-U)^{\ast }(Id_{H}-U)$ \textit{is
compact.}

If $\mathcal{E}$ is the spectral decomposition of $Id_{H}-N$, using
Proposition 2.14, we obtain $E((-\infty ,1-\varepsilon )\cup (1+\varepsilon
,\infty ))=E(g^{-1}((-\infty ,-\varepsilon )\cup (\varepsilon ,\infty )))=%
\mathcal{E(}(-\infty ,-\varepsilon )\cup (\varepsilon ,\infty )\mathcal{)}=%
\mathcal{E(}(-\infty ,-\varepsilon )\cup (\varepsilon ,\infty )\mathcal{)}$,
where $g(x)=1-x$. Since from the above two claims we infer that the operator 
$E(((-\infty ,1-\varepsilon )\cup (1+\varepsilon ,\infty )))=E((-\infty
,1-\varepsilon ))+E((1+\varepsilon ,\infty ))$ has finite rank, we get that $%
\mathcal{E(}(-\infty ,-\varepsilon )\cup (\varepsilon ,\infty )\mathcal{)}$
has finite rank, for every $\varepsilon >0$. Proposition 2.15\ assures us
that $Id_{H}-N$ is compact, i.e. $Id_{H}-(Id_{H}-U)^{\ast
}(Id_{H}-U)=U+U^{\ast }-U^{\ast }U$ is compact.

\textit{Hence} 
\begin{equation*}
S^{[m]}+(S^{[m]})^{\ast }-S^{[m]}(S^{[m]})^{\ast }
\end{equation*}%
\textit{is compact}, \textit{for every }$m\in \mathbb{N}$.

For $m=1$, we get that $S+S^{\ast }-S^{\ast }S$ is compact. Note that, by
Corollary 2.17, $S+S^{\ast }-SS^{\ast }$ is compact and hence $SS^{\ast
}-S^{\ast }S$ is compact (3).

Consequently $S^{\ast }(S^{\ast }S-SS^{\ast })S=(S^{\ast
})^{[2]}S^{[2]}-S^{\ast }SS^{\ast }S$ is compact. (4)

Moreover, for $m=2$, we obtain that $S^{[2]}+(S^{\ast })^{[2]}-(S^{\ast
})^{[2]}S^{[2]}$ is compact. (5)

But%
\begin{equation*}
(S+S^{\ast }-S^{\ast }S)(S+S^{\ast }-S^{\ast }S)=
\end{equation*}%
\begin{equation*}
=(S+S^{\ast })^{[2]}-(S+S^{\ast }-S^{\ast }S)S^{\ast }S-S^{\ast }S(S+S^{\ast
}-S^{\ast }S)-S^{\ast }SS^{\ast }S
\end{equation*}%
is compact.

Since $S+S^{\ast }-S^{\ast }S$ is compact, we infer that%
\begin{equation*}
(S+S^{\ast })^{[2]}-S^{\ast }SS^{\ast }S=
\end{equation*}%
\begin{equation*}
=S^{[2]}+(S^{\ast })^{[2]}+SS^{\ast }+S^{\ast }S-S^{\ast }SS^{\ast }S=
\end{equation*}%
\begin{equation*}
=S^{[2]}+(S^{\ast })^{[2]}-(S^{\ast })^{[2]}S^{[2]}+SS^{\ast }+S^{\ast
}S+(S^{\ast })^{[2]}S^{[2]}-S^{\ast }SS^{\ast }S
\end{equation*}%
is compact. (6)

Then, from (4), (5) and (6), we get that $SS^{\ast }+S^{\ast }S$ is compact.
(7)

From (3) and (7) we deduce that $S^{\ast }S$ is a compact operator and,
using Proposition 2.18, we conclude that $S$ is compact. $\square $

\bigskip

\textbf{4. A necessary condition for the compactness of an operator. }In
this section $X$ is a Banach space. We shall designate by $Id_{X}$ the
function $Id_{X}:X\rightarrow X$, given by $Id_{X}(x)=x$, for all $x\in X$.
If $S$ and $T$ be bounded linear operator from $X$ to $X$ such that $%
\left\Vert S\right\Vert ,\left\Vert T\right\Vert <1$, then $S$ and $T$ are
contractions and $T^{[n]}-Id_{X}$ is invertible, for each $n\in \mathbb{N}$.
For $w\in X$, we consider the IFS $S_{w}=(X,f_{1},f_{2})$, where $%
f_{1},f_{2}:X\rightarrow X$ are given by $f_{1}(x)=S(x)$ and $%
f_{2}(x)=T(x)+w $, for all $x\in X$.

\bigskip

THEOREM\ 4.1. \textit{In the above mentioned setting, if the operator }$S$ 
\textit{is compact, then} \textit{there exists a family }$(K_{n})_{n\in 
\mathbb{N}}$ \textit{of compact subsets of }$X$\textit{\ such that }$A_{%
\mathcal{S}_{w}}$\textit{\ is\ not connected, for all }$w\in H-\underset{%
n\in \mathbb{N}}{\cup }K_{n}$\textit{.}

\bigskip

$\mathtt{Proof}$. The proof given in Theorem 5, from [15], applies with
little change. More precisely let $C_{0}$ be the compact set $\overline{%
S(B(0,1))}$. Let $X^{^{\prime }},X_{1},X_{2},...,X_{n},...$ be given by 
\begin{equation*}
X^{^{\prime }}=S(X)=\underset{k\in \mathbb{N}}{\cup }kC_{0}
\end{equation*}%
and 
\begin{equation*}
X_{n}=(T-Id_{X})(T^{[n]}-Id_{X})^{-1}(X^{^{\prime }}-T^{[n]}(X^{^{\prime }}))%
\text{,}
\end{equation*}%
for each $n\in \mathbb{N}$. We have 
\begin{equation*}
X_{n}=(T-Id_{X})(T^{[n]}-Id_{X})^{-1}(\underset{k\in \mathbb{N}}{\cup }%
kC_{0}-T^{[n]}(\underset{l\in \mathbb{N}}{\cup }lC_{0}))=
\end{equation*}%
\begin{equation*}
=(T-Id_{X})(T^{[n]}-Id_{X})^{-1}(\underset{k\in \mathbb{N}}{\cup }kC_{0}-%
\underset{l\in \mathbb{N}}{\cup }lT^{[n]}(C_{0}))=
\end{equation*}%
\begin{equation*}
=(T-Id_{X})(T^{[n]}-Id_{X})^{-1}(\underset{k,l\in \mathbb{N}}{\cup }%
(kC_{0}-lT^{[n]}(C_{0}))\text{,}
\end{equation*}%
for each $n\in \mathbb{N}$ and since $kC_{0}-lT^{[n]}(C_{0})$ is compact for
all $k,l\in \mathbb{N}$, we infer that $X_{n}$ is a countable union of
compact subsets of $X$. Therefore there exists a family\textit{\ }$%
(K_{n})_{n\in \mathbb{N}}$ of compact subsets of $X$ such that $\underset{%
n\in \mathbb{N}}{\cup }X_{n}=\underset{n\in \mathbb{N}}{\cup }K_{n}$. The
rest of the proof of the Theorem mentioned above does not require any
modification.

Hence $A_{\mathcal{S}_{w}}$\textit{\ }is disconnected, for each $w\in
X\smallsetminus \underset{n\in \mathbb{N}}{\cup }X_{n}=X\smallsetminus 
\underset{n\in \mathbb{N}}{\cup }K_{n}$. $\square $

\bigskip

REMARK 4.2. If $X$ is infinite dimensional, then $W\overset{not}{=}%
X\smallsetminus \underset{n\in \mathbb{N}}{\cup }X_{n}=X\smallsetminus 
\underset{n\in \mathbb{N}}{\cup }K_{n}$ is dense in $X$.

\bigskip

$\mathtt{Proof}$. Indeed, let us note that $K_{n}$ is a closed set. Moreover 
$\overset{\circ }{K_{n}}=\varnothing $ since if this is not the case, then
the closure of the unit ball of the infinite-dimensional space $X$ is
compact which is a contradiction. Consequently $X_{n}$ is nowhere dense, for
each $n\in \mathbb{N}$, and therefore $W$ is dense in $X$.

\bigskip

\begin{center}
\textbf{References}
\end{center}

\bigskip

[1] D. Baki\'{c} and B. Gulja\v{s}, \textit{Which operators approximately
annihilate orthonormal bases?}, Acta Sci. Math. (Szeged) 64 (1998), No.3-4,
601-607.

[2] M.F. Barnsley, \textit{Fractals everywhere}, Academic Press
Professional, Boston, 1993.

[3] J. J. M. Chadwick and A.W. Wickstead, \textit{A quotient of ultrapowers
of Banach spaces and semi-Fredholm operators}, Bull. London Math. Soc. 9
(1977), 321-325.

[4] J. B. Conway, \textit{A course in functional analysis}, Springer-Verlag,
New York, Berlin, Heidelberg, Tokyo, 1985.

[5] R. Cristescu, \textit{No\c{t}iuni de Analiz\u{a} Func\c{t}ional\u{a}
Liniar\u{a}}, Editura Academiei Rom\^{a}ne, Bucure\c{s}ti, 1998.

[6] D. Dumitru and A. Mihail, \textit{A sufficient condition for the
connectedness of the attractors of an infinite iterated function systems},
An. \c{S}tiin\c{t}. Univ. Al. I. Cuza Ia\c{s}i. Mat. (N.S.), LV (2009), f.1,
87-94.

[7] K.J. Falconer, \textit{Fractal geometry, Mathematical foundations and
applications}, John Wiley \& Sons, Ltd., Chichester, 1990.

[8] P. A. Fillmore and J.P. Williams, \textit{On operator ranges}, Adv.
Math. 7 (1971), 254-281.

[9] J.E. Hutchinson, \textit{Fractals and self-similarity}, Indiana Univ.
Math. J. 30 (1981), 713-747.

[10] K. Imazeki, \textit{Characterizations of compact operators and
semi-Fredholm operators}, TRU Math. 16 (1980), No.2, 1-8.

[11] J. M. Isidro, \textit{Structural characterization of compact operators}%
, in Current topics in operator algebras, (Nara, 1990), 114-129, World Sci.
Publ., River Edge, NJ, 1991.

[12] J. Kigami, \textit{Analysis on Fractals,} Cambridge University Press,
2001.

[13] R. Miculescu and A. Mihail, \textit{Lipscomb's space }$\omega ^{A}$%
\textit{\ is the attractor of an infinite IFS containing affine
transformations of }$l^{2}(A)$, Proc. Amer. Math. Soc. 136\ (2008), No. 2,
587-592.

[14] A. Mihail, \textit{On the connectivity of the attractors of iterated
multifunction systems}, Real Anal. Exchange, 34 (2009), No. 1, 195-206.

[15] A. Mihail and R. Miculescu, \textit{On a family of IFSs whose
attractors are not connected}, to appear in J. Math. Anal. Appl., DOI
information: 10.1016/j.jmaa.2010.10.039.

[16] K. Muroi and K. Tamaki, \textit{On Ringrose's characterization of
compact operators}, Math. Japon. 19 (1974), 259-261.

[17] J. R. Ringrose, \textit{Compact non-self-adjoint operators}, Van
Nostrand Reinhold Company, London, 1971.

[18] \'{A}. Rod\'{e}s-Us\'{a}n, \textit{A characterization of compact
operators}, in Proceedings of the tenth Spanish-Portuguese on mathematics,
III (Murcia, 1985), 178-182, Univ. Murcia, Murcia, 1985.

[19] W. Rudin, \textit{Functional analysis}, 2nd ed., International Series
in Pure and Applied Mathematics. New York, NY: McGraw-Hill, 1991.

\bigskip

Department of Mathematics

Faculty of Mathematics and Informatics

University of Bucharest

Academiei Street, no. 14

010014 Bucharest, Romania

E-mail: mihail\_alex@yahoo.com

\qquad\ \ \ \ \ miculesc@yahoo.com

\end{document}